\documentclass[12pt]{amsart}
\usepackage{lineno}

\usepackage{amssymb, amscd}
\usepackage{mathptmx}
\usepackage{lmodern}

\usepackage{xcolor}
\definecolor{cornsilk}{rgb}{1.0, 0.97, 0.86}

\usepackage{soul}
\usepackage{mathrsfs}
\numberwithin{equation}{section}

\DeclareSymbolFont{SY}{U}{psy}{m}{n}
\DeclareMathSymbol{\emptyset}{\mathord}{SY}{'306}

\newcommand{\overbar}[1]{\mkern 1.5mu\overline{\mkern-1.5mu#1\mkern-1.5mu}\mkern 1.5mu}

\usepackage{todonotes}

\usepackage{tikz-cd}

\theoremstyle{plain}

\newtheorem{thm}{Theorem}[section]

\theoremstyle{definition}

\newtheorem{rem}[thm]{Remark}

\DeclareMathOperator{\ad}{ad}

\newcounter{defcounter}
\makeatletter
\def\moverlay{\mathpalette\mov@rlay}
\def\mov@rlay#1#2{\leavevmode\vtop{%
   \baselineskip\z@skip \lineskiplimit-\maxdimen
   \ialign{\hfil$\m@th#1##$\hfil\cr#2\crcr}}}
\newcommand{\charfusion}[3][\mathord]{
    #1{\ifx#1\mathop\vphantom{#2}\fi
        \mathpalette\mov@rlay{#2\cr#3}
      }
    \ifx#1\mathop\expandafter\displaylimits\fi}
\makeatother

\input xy
\xyoption{all}
\usepackage{fancyhdr}
\makeatletter

\title
{A simplified approach to the holomorphic discrete series}

\author[A. Kor\'{a}nyi]{Adam Kor\'{a}nyi} 

\address[A. Kor\'{a}nyi]{Lehman College\\
Bronx, NY 10468}

\email[A. Kor\'{a}nyi]{adam.koranyi@lehman.cuny.edu}

\begin{document}

\begin{abstract}
A simple new proof of the Harish-Chandra condition, preceded by an expository part on Hermitian symmetric spaces, holomorphic induction, and on some analytic tools.
\end{abstract}
\maketitle 
\setcounter{section}{-1}
\section{Introduction} 
The only new thing in this article is the short proof in $\S\, 6$  of the ``Harish-Chandra condition" for the existence of the holomorphic discrete series. A simpler proof than the original one \cite{HC1}, \cite{HC2}, was found by H. Rossi and M. Vergne in 1975 \cite{RV} as a special case of a more general theory, but it is still considerably more complicated than the one given here. I noticed the possibility of such a short proof a few years ago while studying an old article of Wallach \cite{W2} and sent a brief report on it to arXiv \cite{K3}.

With the present article I have a double purpose. One is to make the simple proof easily readable and accessible. The other is to write a fairly detailed expository paper, with at least outlines or indications of proofs, about the bounded symmetric domains and the holomorphic discrete series from the analyst's point of view. This includes the basic facts about holomorphic induction, the Harish-Chandra imbedding theorem and a section on reproducing kernel spaces. I will not go beyond proving the existence condition, so exact values of positive constants will often not be computed; several of these can be found in the second  section of \cite{K2}. Analytic continuation of the representations, formal dimension, characters will not be included either, except for some brief remarks in $\S\, 7$.  

What is assumed as known is first of all the basic theory of semisimple Lie algebras, their real forms, roots and weights. About this there are  a number of good specialized books; most of what we need is also in \cite{H1} and \cite{W1}. Furthermore, the basic facts about the globally symmetric space $G / K$, the decomposition $G=K A K$, etc. are also taken for granted. These matters (along with much else) are mostly to be found in \cite{H1} and \cite{W1}. I also have an expository (partly just heuristic) introduction to the subject \cite{K1}; it has a slight overlap with the present article as well. In the present article, in $\S\, 1,$ the discussion begins with the case where the symmetric space has an invariant complex structure.

\section{Algebraic preliminaries.} Our starting point is a simple real Lie algebra $\mathfrak g$ in whose Cartan decomposition $\mathfrak g = \mathfrak k + \mathfrak p$, the subalgebra $\mathfrak k$ is not semisimple. (This is equivalent to the existence of a $G$-invariant complex structure on the symmetric space.) Then $\mathfrak k=\mathfrak k_{\text{ss}} \oplus \mathfrak z$ with $\mathfrak k_{\text{ss}}=[\mathfrak k, \mathfrak k]$ semisimple and $\mathfrak z$ the (necessarily one-dimensional) center. $\mathfrak g^\mathbb C$ is the complexification of $\mathfrak g$, and $\mathfrak g_U=\mathfrak k+i \mathfrak p$ is a compact real form. Writing $B$ for the Killing form and $\tau$ for conjugation with respect to $\mathfrak g_U$, $B_\tau(x, y)=-B(x, \tau y)$ is a Hermitian inner  product on $\mathfrak g^{\mathbb C}$. We fix a maximal Abelian subalgebra $\mathfrak h$ in $\mathfrak k$; we have $\mathfrak h= \mathfrak h_{\text{ss}} \oplus \mathfrak z$ with $\mathfrak h_{\text{ss}} \subset \mathfrak{k}_{\text{ss}}$.

The set $\Delta$ of $\mathfrak h^\mathbb C$-roots $\alpha$ have corresponding root spaces $\mathfrak g_\alpha$; each of these is either in $\mathfrak{k}^\mathbb C$ or in $\mathfrak p^\mathbb C$ (by $[\mathfrak h, \mathfrak k] \subset \mathfrak k$, etc.). Accordingly, we have the sets $\Delta_c$ and $\Delta_n$ of ``compact" resp.\ ``non-compact" roots. We choose an ordering, i.e. a set $\Delta^{+}$of positive roots, in such a way that on any fixed generator of $i \mathfrak z$ the roots in $\Delta_n^+ =  \Delta^+\cap \Delta_n$ are all positive or all negative.  We have also the sets $\Delta^{-}, \Delta_c^{+}, \Delta_c^{-}$, etc.

We fix a 
basis i.e., for each $\alpha \in \Delta^{+}$, elements $h_\alpha \in i \mathfrak h, e_{ \pm \alpha} \in \mathfrak g_{ \pm \alpha}$ such that
$$
\left[h_\alpha, e_\alpha\right]=2 e_\alpha, \quad \tau e_\alpha=-e_{-\alpha}, \quad\left[e_\alpha, e_{-\alpha}\right]=h_\alpha. 
$$
(To make the conventions clear: For an element $\varphi$ of the real dual of $i\mathfrak{h}$, we have $\varphi(h_\alpha)= 2(\varphi \mid\alpha)/(\alpha\mid\alpha)$ for any Weyl group invariant inner product, e.g. for the restriction of $B_\tau$.)
The direct sums of the root spaces for $\Delta_n^{+}, \Delta_n^{-}$ give Abelian subspaces $\mathfrak p^{+}, \mathfrak p^{-}$ of $\mathfrak g^\mathbb C$ such that $\mathfrak p^c=\mathfrak p^{+} \oplus \mathfrak p^{-}$. There is an element $\hat{z}$ in $\mathfrak{z}$ such that $\operatorname{ad}(\hat{z})= \pm i I$ on $\mathfrak p^{ \pm}$. Then $J=\operatorname{ad}(\hat{z})$ is a complex structure on $\mathfrak p$, and $j=\frac{1}{2}(I-i J)$ is an $\operatorname{ad}(\mathfrak k)$-equivariant isomorphism $\mathfrak p \to \mathfrak p^+$.

There is a crucially important way to construct a maximal Abelian subalgebra $\mathfrak a$ of $\mathfrak p$ (\cite[Sec.6]{HC1}, \cite[VIII., Prop 7.4]{H1}).  
Starting with $\gamma_r$, the highest root in $\Delta_n^+$ (more convenient here than the original order starting with $\gamma_1$, the lowest) one constructs a set of roots $\gamma_1, \ldots, \gamma_r$ in $\Delta_n^{+}$ such that $\gamma_j \pm \gamma_k$ is never a root. Using the abbreviations
$$
h_j=h_{\gamma_j}, \quad e_{ \pm j}=e_{ \pm \gamma_j}, \quad(1 \leq j \leq r)
$$
the subalgebra
\begin{equation}\label{eqn:1.1}
\mathfrak a =\sum_{j=1}^{r} \mathbb{R}\left(e_j+e_{-j}\right)
\end{equation}
is maximal in $\mathfrak p$. ($r$ is the rank of $G / K$). The key fact here is that the centralizer of $e_r+e_{-r}$ in $\mathfrak{p}_\mathbb{C}$ is the direct sum of $\mathbb{C}(e_r+e_{-r})$ and the space $\mathfrak p_{r-1}$, defined as the sum of all $ \mathfrak{g}_\gamma + \mathfrak {g}_{-\gamma} \: (\gamma \in \Delta_n^+)$  such that $\gamma \pm \gamma_{r}\: $ are not roots. Then one can naturally define a subalgebra $\mathfrak{g}_{r-1}^\mathbb{C}$ (and $\mathfrak g_{r-1})$ which gives again a Hermitian symmetric space (of rank $r-1$)  whose $\mathfrak{p}$-part is $\mathfrak{p}_{r-1}.$ One can repeat the process, getting a chain of subalgebras $\mathfrak{g}_j \: (1 \leq j \leq r)$
of rank $j$. (These even give a chain of Hermitian symmetric spaces, each geodesically imbedded in the preceding one.)

$\mathfrak h$  splits into $\mathfrak h^{-}=\sum_1^{r} \mathbb{R} h_j$ and its  orthocomplement $\mathfrak h^{+}$. There is an inner automorphism (a ``Cayley transformation" $\operatorname{Ad}(\operatorname{exp}\frac{\pi}{4}\sum(e_j - e_{-j})$) in $\mathfrak g^\mathbb{C}$ leaving $i\mathfrak {h}^+$ pointwise fixed and mapping $\mathfrak a$ to $i \mathfrak h^-$. Hence the system of $\mathfrak a$-roots is isomorphic with the restrictions of $\mathfrak h^\mathbb C$-roots to $\mathfrak i \mathfrak h^-$; briefly the ``restricted roots". (Recall that $i\mathfrak h$ is the space where the roots are real; the root system is a subset of the real dual space of $i\mathfrak h^-$.)

Up to here everything is in \cite{H1}. It was not yet known at the time but became known later that the $\mathfrak a$-roots always form an abstract root system. This 
together with some of Harish-Chandra's original lemmas made it easy to see that in our case the restricted root system is of type $B_r$ or $B C_r$ \cite[V., $\S\, 4$]{H3}. The restricted roots coming from $\Delta_n^+$ are the $\gamma_j$ $(1 \leq j \leq r)$, each with multiplicity 1; $\frac{1}{2}\left(\gamma_j+\gamma_k\right)\,\,(1 \leq k<j \leq r)$ all with the same multiplicity $a$; and $\frac{1}{2} \gamma_j$ with multiplicity $b$ (which may be $0$). The restricted roots from $\Delta_c^{+}$ are $\frac{1}{2}\left(\gamma_j-\gamma_k\right)$ again with multiplicity $a$, and $\frac{1}{2} \gamma_j$ with multiplicity $b$. 

The Weyl group consists of the signed permutations of the roots $\gamma_j$. 

An important invariant of the space is
\begin{equation} \label{eqn:1.2}
p=(r-1) a + b + 2
\end{equation}
Its significance is due to the fact that if we make the usual definitions
\[
\rho=\frac{1}{2} \sum_{\alpha \in \Delta^{+}} \alpha, \quad \rho_n=\frac{1}{2} \sum_{\alpha \in \Delta^{+}_n} \alpha,
\]
it follows easily from what was said above that
\begin{align}
\rho\left(h_{r}\right)&=p-1 \label{eqn:1.3}\\
2 \rho_n\left(h_j\right)&=p \quad(1 \leq j \leq r). \label{eqn:1.4}
\end{align}
\section{The Harish-Chandra imbedding} We denote by $G^{\mathbb{C}}$ the simply connected group with Lie algebra $\mathfrak g^{\mathbb{C}}$, by $G_0$, $K_0$, $Z_0$, $K_0^{\mathbb{C}}$, $G_U$, $P^{+}$, $P^{-}$, $A$ its analytic subgroups for $\mathfrak g, \mathfrak k, \mathfrak z, \mathfrak k^\mathbb{C}, \mathfrak {g}_U, \mathfrak p^{+}, \mathfrak p^-, \mathfrak a$. (The subscripts $0$ are there on the first four subgroups because we reserve the notations $G, K, Z, K^{\mathbb C}$ for their simply connected coverings.) We note that $P^+, P^{-}$ are closed Abelian subgroups, isomorphic with $\mathfrak p^+, \mathfrak p^-$ under the exponential map. $K_0^\mathbb C P^+$ and $K_0^\mathbb C  P^-$ are also closed subgroups. $\mathcal M=G_0 / K_0$ is the symmetric space we want to realize as a bounded domain.

It is a crucial, even if not very difficult, fact that $P^{+} K^{\mathbb{C}} P^{-}$ is open in $G^{\mathbb{C}}$ and the product map, from $P^{+} \times K^{\mathbb{C}} \times P^{-}$ is a holomorphic diffeomorphism (\cite[VIII., Lemma 7.9]{H1}).



We consider the homogeneous complex analytic manifold $\mathcal M^*=G^\mathbb C / K_0^{\mathbb C} P^{-}$, and we take in it the orbits of the base point, successively, under $G_U$, $P^+$, and $G_0$. This gives imbeddings of $G_U/ K_0$, $P^+$, and $\mathcal M=G_0/K_0$ into $\mathcal M^*$ as open complex submanifolds (because $G_U \cap K_0^\mathbb C P^{-}=K_0$, etc.). For $G_U/ K_0$ the image is also compact besides being open, so $G_U /K_0 \cong \mathcal M^*$. So we have an imbedding of $\mathcal M$ into its compact dual $\mathcal M^*$; this is known as the Borel imbedding.

Our focus will now be on the imbedding of $ P^+$ into $\mathcal M^*$ which also gives an imbedding (via the exponential) of the complex vector pace $\mathfrak p^+$. We will show that the image of $\mathcal M$ is contained in the image of $\mathfrak p^+$, thus realizing $\mathcal M$ as a domain in $\mathfrak p^{+}$ with $G_0$ acting on it by holomorphic mappings.

To formalize this, for $z$ in $\mathfrak p^{+}$ and $g \in G^\mathbb C$ we define $g(z)$, or briefly $gz$, provided $g \exp z$ is in $P^{+} K_0^\mathbb C P^{-}$, by writing the $P^+K_0^\mathbb C P^-$  decomposition as
\begin{equation} \label{eqn:2.1}
g \exp z=(\exp g z) J(g, z) p^-(g, z). 
\end{equation}
The first factor is the definition of $g z$.

The second factor is Satake's ``canonical automorphy factor'' \cite[Ch. II, $\S\, 5$]{S}, and will play a fundamental role. When $gz$, $g_1 z$, $gg_1 z$ are defined, it is easy to verify that
\begin{equation} \label{eqn:2.2}
J\left(g g_1, z\right)=J\left(g, g_1 z\right) J(g_1, z),
\end{equation}
so $J(g, z)$ is a $K_0^{\mathbb{C}}$ valued cocycle (or multiplier, or automorphy factor). Furthermore,
\begin{equation} \label{eqn:2.3a}
J(k,z) = k\end{equation}
for all $k\in K_0^\mathbb C, z\in \mathfrak p^+$

Returning to the action on $\mathfrak p^{+}$, it is immediately clear that elements of $P^+$ act by translation and that for $g=k \in K_0^{\mathbb{C}}$ we have
\begin{equation} \label{eqn:2.3}
k z=\operatorname{Ad}(k) z,
\end{equation}
a unitary transformation.

Now we are going to use the decomposition $G_0=K_0 A K_0$. We parametrize  $A$ by $t=(t_1, \ldots, t_r) \in \mathbb{R}^{r}$ as
\begin{equation} \label{eqn:2.4}
a(t)=\exp \sum_1^r t_j\left(e_j+e_{-j}\right).
\end{equation}

For fixed $j$ we compute in (any) Lie group with Lie algebra 
$\mathbb C e_j + \mathbb C h_j + \mathbb C e_{-j}$ the $P^+ K_0^\mathbb C P^-$ decomposition:
\begin{equation} \label{eqn:2.5}
\operatorname{exp} t_j\left(e_j+e_{-j}\right)=\left(\exp \left(\tanh t_j\right ) e_j\right)\left(\exp \left(-\log \cosh t_j\right) h_j\right)\left(\exp \left(\tanh t_j\right) e_{-j}\right).
\end{equation} 

For different $j$-s the $3$-dimensional algebras completely commute with each other, since $\gamma_j \pm \gamma_k$ is never a root.

Consequently, one can take the product of the identities \eqref{eqn:2.5} over all $j$. This product then gives an explicit expression of \eqref{eqn:2.1} for $g=a(t)$, $z=0$, and yields
\begin{equation}\label{eqn:2.6}
a(t) \cdot 0=\sum_i^p\left(\tanh t_j\right) e_j \\
\end{equation}
\begin{equation}\label{eqn:2.7}
J(a(t), 0)=\exp  \big (-\sum_{1}^r  \log \cosh t_j\big ) h_j.
\end{equation}
Now the multiplier identity and \eqref{eqn:2.3a} immediately give 
\begin{equation} \label{eqn:2.7a}
J(k a(t),0) = k J(a(t),0)
\end{equation}

As we see from \eqref{eqn:2.6}, the full orbit of 0 under $A$ is an $r$-dimensional real cube (with corners $(\pm e_1, \ldots , \pm e_r ))$. The orbit of $G_0$ is this cube rotated around in $\mathfrak p^+$ by the group $\operatorname{Ad}(K_0)$ of unitary transformations. This is $D$, the space $\mathcal M$ holomorphically imbedded in $\mathfrak p^{+}$ as a bounded domain. $G_0$ acts on it by holomorphic self -maps, there is a smooth invariant Riemannian metric and a corresponding $G_0$-invariant measure $\nu$ on $D$. There is also, as on the ambient $\mathfrak p^{+}$, the Euclidean measure $\lambda$ on $D$; its relation to $\nu$ is computed in $\S\, 5$.

Since we will mainly work with $D$, it is useful to also have and use the notation $x_j= \tanh t_j$ $(1 \leq j \leq r)$. With that, \eqref{eqn:2.6} and \eqref{eqn:2.7} take the form 
\begin{equation}\label{eqn:2.8}
 a(t) \cdot 0=\sum_1^r x_j e_j \\
\end{equation}
\begin{equation}\label{eqn:2.9}
 J(a(t), 0)=\exp \left(\frac{1}{2} \sum_1^r \log (1-x_j^2) h_j\right)
\end{equation}
(by $\cosh ^ {-2} t= 1 -\tanh ^2 t$ ).

The $K_0 A K_0$ decomposition also leads to a formula of integration in generalized polar coordinates. As $\lambda$ is $K_0$-invariant, and the $K_0$-orbits are orthogonal to $A\cdot 0$ in $\mathfrak p^{+}$, it is not too difficult to prove that
\begin{equation}\label{eqn:2.10}
 \int_D f d \lambda =c \int_{K_0 /M_0} d k_M \int_{0\leq x_1 \leq \cdots \leq x_r  < 1} d x f\left ( \operatorname{Ad}(k) \sum x_j e_j\right) P(x),
 \end{equation}
where $c>0$, $M_0$ is the centralizer of $\mathfrak a$  in $K_0$, $d k_M$ is the normalized $K_0$-invariant measure on $K_0 /M_0$, the $x$-integral is on a fundamental domain of the Weyl group, and
$$
P(x)=\prod_{j=1}^r x_j^{2 b+1} \prod_{1 \leq j<k \leq r}\left(x_k^2-x_j^2\right)^a
$$
is the product of the restricted roots. (This is Theorem 5.7 in \cite[Ch I]{H2} transferred by $j$ from $\mathfrak p$ to $\mathfrak p^+$; cf. also \cite[Sec. 7.8]{W1}.)

\section{The holomorphic discrete series} We want representations of the simply connected covering group $\zeta: G \rightarrow G_0$, so we have to look at $D$ as $G / K$ ($K$ the analytic subgroup for $\mathfrak k$). The action of G on $D$ is then given by $g z=\zeta(g) z$. We also need $J(g, z)$ defined by \eqref{eqn:2.1} to be continuous on $G \times D$ and still satisfy \eqref{eqn:2.2}. For this we have to lift $J(g, z)$ to take values in the simply connected covering group $\xi: K^{\mathbb{C}} \rightarrow K_0^{\mathbb{C}}$. In the following we use $J(g, z)$ in this sense, without changing the notation. (For more on this see Remark 1 in $\S\,7$.)

We define the representations in the holomorphic discrete series as representations holomorphically induced by irreducible unitary representations $(\pi, V)$ of $K$. This means the restriction of the unitarily induced (Mackey) representation $\operatorname{Ind}_K^G(\pi)$ to the holomorphic subspace with respect to a $G$-invariant complex structure on the homogeneous vector bundle $G\times_\pi V$. In our case, there is only one such complex structure, because the Lie algebra $\mathfrak k^\mathbb C+\mathfrak p^{-}$ which characterizes the complex structure of $\mathcal M \simeq D$ has only one representation extending $\pi |_{\mathfrak k^\mathbb C}$, namely 0 on $\mathfrak p^{-}$.
(For clear descriptions of these matters see \cite[Ch. I]{V}.)

There is a simple way to realize $\operatorname{Ind}_K^G(\pi)$ as a multiplier representation. We fix $\pi$ and do not indicate in the notation the dependence on it of the subsequent definitions. We set
$$
m_g(z)=m(g, z)=\pi(J(g, z)).
$$
$m$ is then a $\operatorname{GL}(V)$-valued multiplier, i.e. satisfies a cocycle identity like \eqref{eqn:2.2} on $G \times D$; this implies the properties $m(e, z) \equiv I$, $m\left(g, g^{-1} z\right)=m\left(g^{-1}, z\right)^{-1}$. We also have
$$
m(k, z)=\pi(k) \quad\left(k \in K^{\mathbb{C}}, z \in D\right). 
$$

$m$ gives us the ``multiplier representation" $U$ of $G$ defined as $U_g f=\left(m_g f\right) \circ g^{-1}$, or, a way of writing it out, 

\begin{equation}\label{eqn:3.1}
\big(U_g f\big)(z)=m\left(g^{-1}, z\right)^{-1} f\left(g^{-1} z\right).
\end{equation}
It is unitary on the Hilbert space $L^2(D, \pi, Q, \nu)$ of $V$-valued functions with the norm
\begin{equation}\label{eqn:3.2}
\|f\|^2=\int_D\langle Q(z) f(z), f(z)\rangle_V d\nu(z),
\end{equation}
where $Q(z)$ is a positive definite element of $\operatorname{GL}(V)$, continuous in $z$ and well defined by
\begin{equation}\label{eqn:3.3}
Q(g \cdot 0)=m(g, 0)^{*-1} m(g, 0)^{-1}.
\end{equation}
As before, $\nu$ is the invariant measure.

We note that $Q(0)=I$, and for all $g \in G, z \in D$,
\begin{equation}\label{eqn:3.4}
Q(g z)={m(g, z)^*}^{-1} Q(z) m(g, z)^{-1}
\end{equation}
as one sees writing $z=g_1.0$ and applying the multiplier identity. In particular,
\begin{equation}\label{eqn:3.5}
Q(k z)=\pi(k) Q(z) \pi(k)^{-1} \quad(z \in D, k \in K).
\end{equation}
Note also that $C_c^{\infty}(D, V)$ is dense in $L^2(D, \pi, Q, \nu)$, by the usual arguments.

We want to see that $U$ is indeed a realization of $\operatorname{Ind}_K^G(\pi)$. Here, with everything simply connected and smooth, this is almost trivial: One usually defines $\operatorname{Ind}_K^G(\pi)$  by taking the space $C_0^{\infty}(G, V)^\pi$ of functions $F: G \rightarrow V$ whose projection to $D$ is $C^{\infty}$ with compact support, and satisfy
$$
F(g k)=\pi(k)^{-1} F(g)
$$
(i.e.  lifts to $G$ of sections of $G \times_\pi V$), and put on them the norm
$$
\|F\|^2=\int_D |F |_V^2(z) d\nu(z),
$$
where we write $|F|^2(g\cdot  0)=|F(g \cdot 0)|_V^2$, which depends only on $g\cdot 0$, the $K$-coset of $g$. $\operatorname{Ind}_K^G(\pi)$ is left translation by $G$ on the completed Hilbert space.
Now the relation
$$
F(g)=m(g, 0)^{-1} f(g \cdot 0)
$$
gives a unitary equivalence of representations, proving our assertion.

Getting to \textit{holomorphic} induction: Let $\mathcal H$ be the set of holomorphic functions in $L^2(D, \pi, Q, \nu)$ (in the ordinary sense), Since $m(g, z)$ is holomorphic in $z$, the multiplier representation $U$ preserves $\mathcal H$. In $\S\,4$. we show (from an independent start) that $\mathcal H$ is actually a closed subspace, so we can conclude: The holomorphic discrete series is the restriction of $U$ to $\mathcal H$.   

The only question is for what $\pi$ is the space $\mathcal H$  non-zero; it is answered by the Harish-Chandra criterion, in $\S\, 6$.
\section{Reproducing kernel spaces} For scalar-valued functions all that follows is very well known. We describe the $V$-valued case ($V$ a finite dimensional Hilbert space), which is an easy generalization. Let $E$ be a set, and let $H$ be a Hilbert space of functions (not a.e classes!) $f: E \rightarrow V.$ By a kernel we mean a function $\mathcal K: E \times E \rightarrow  \operatorname{GL}(V)$, and use the notations
\begin{align*}
\mathcal K_w(z) = \mathcal K(z, w), & \quad z, w \in E; \\
\left(\mathcal K_w v\right)(z)=\mathcal K(z, w) v, & \quad v \in V.
\end{align*} 

$\mathcal K$ is a \textit{reproducing kernel} for $H$ if  $\mathcal K_w v \in H$ for all $w \in E$, $v \in V$, and
\begin{equation}\label{eqn:4.1}
\left(f \mid \mathcal K_w v\right)=\langle f(w), v\rangle_V
\end{equation}
for all $f \in H$. The Riesz representation theorem implies that such a kernel exists (and is then unique) if and only if the evaluation functionals  $\operatorname{ev}_w: H \rightarrow V$, given by $\operatorname{ev}_w(f)=f(w)$, are continuous for all fixed $w$. Note that the vectors $\mathcal K_w v$ span $H$; also that $\mathcal K_w=\operatorname{ev}_w^*$  and $\mathcal K\left(z, w\right)=\operatorname{ev}_z \circ \operatorname{ev}_w^*$. For easy reference we write down the properties
\begin{equation}\label{eqn:4.2}
\left(\mathcal K_w v \mid \mathcal K_z v^{\prime}\right)=\big\langle\mathcal{K}(z, w) v, v^{\prime}\big \rangle_V,
\end{equation}
\begin{equation}\label{eqn:4.3}
\mathcal K(z, w)=\mathcal K(w, z)^*.
\end{equation}
(here $*$ means the adjoint in $V$).

To apply this to $\mathcal H$, the holomorphics in $L^2(D, \pi, Q, \nu)$, we must show that $\mathcal H$ is a closed subspace and $\operatorname{ev}_w$ is continuous on it. We start with the latter, but do it right away uniformly for $w$ in an (arbitrarily) fixed compact set $S\subset D$. Then there exists $r>0$ such that the polydisc $P_r(w)$ centered at $w \in S$ is still contained in some compact $S^{\prime}$. The (Euclidean) volume of $P_r(w)$ is independent of $w$, so it follows from the mean value theorem for holomorphic functions that with some $c_1>0$,
$$
|f(w)|_V ^2 \leq c_1 \int_{P_r(w)}|f|_V^2 d \lambda
$$
for $w\in S$. By compactness and continuity there is another $c_2 \geq 0$ such that $Q(w) \succ c_2 I$ in the sense of positivity of operators, and $c_3>0$ such that $d\nu \geq c_3 d \lambda$ on $S^{\prime}$. So finally
\[
|f(w)|^2 \leq c_1 \int_{s^{\prime}}|f|_V^2 d \lambda \leq c_1 c_2^{-1} c_3^{-1} \int_{s^{\prime}}\langle Q(z) f(z),f(z)\rangle_V d\nu(z) \leq c_1 c_2^{-1} c_3^{-1} \| f \|^{2}\]
This proves the continuity of evaluations, and also proves that convergence in norm implies uniform convergence on compact sets. So $\mathcal H$  is closed, and (and if non-zero) has a reproducing kernel $\mathcal K$.

To find out more about $\mathcal K$  one way is to use the unitarity of $U$, which gives
$$
\left(U_g \mathcal K_w v \mid \mathcal K_{g z} v^{\prime}\right)=\left(\mathcal K_w v \mid U_g^{-1} \mathcal K_{g z} v^{\prime}\right). 
$$
Then the definition of $U$ and the reproducing property \eqref{eqn:4.1} give after some rearrangement 
\begin{equation}\label{eqn:4.4}
\mathcal K(g z , gw) = m(g,z) \mathcal K(z,w) m(g,w)^*\quad (g\in G; z,w \in D). \end{equation}  
 Since $\mathcal K(z, w)$ is holomorphic in $z$ and antiholomorphic in $w$, it is completely determined by $\mathcal{K}(z, z)$ via the power series
\begin{equation}\label{eqn:4.5}
 \mathcal K(z, z)=\sum c_{\alpha \beta} z^\alpha \overbar{z}^\beta
\end{equation} 
($\alpha$ and $\beta$ the usual multi-indices). By \eqref{eqn:4.3} we have $c_{\beta \alpha}=c_{\alpha \beta}^*$. Setting $w=0$ in the power series of $\mathcal K(z,w)$ we see that
\begin{equation}\label{eqn:4.6} 
\mathcal K(z, 0)= c I \quad (\forall z \in D)
\end{equation}
with some constant $c>0$.

Comparing \eqref{eqn:4.4} with \eqref{eqn:3.4} we have
\begin{equation}\label{eqn:4.7} 
\mathcal K(z, z)=c Q(z)^{-1}. 
\end{equation}

We should note that there is a conceptually simpler way (but involving slightly more computation) to get here. Once we know that $\mathcal K$ exists, we can use the ``polarization" process \eqref{eqn:4.5} on $Q(z)^{-1}$, and verify that it is (up to constant) the reproducing kernel. There is even a third way: To check from \eqref{eqn:3.3} that $Q(z)=c^{\prime} \pi\left(\mathbb{K}(z, z)^{-1}\right)$ with the``canonical kernel function" $\mathbb K$ of Satake, and use the properties of $\mathbb{K}$ proved in \cite[Ch. II, $\S\,5$]{S}.

One consequence of \eqref{eqn:4.7}, very important for us, is the non-trivial fact that if $\mathcal H \neq 0$, then for all $v \in V$, the constant functions $c^{-1}\left(\mathcal K_0 v  \right)(z) \equiv v$ are in $\mathcal H$.

Another consequence of our considerations is the irreducibility of $U$ on $\mathcal H$ (if $\neq 0$). Indeed, suppose $\mathcal H^\prime \subset \mathcal H$ is a closed invariant subspace of $U$, which necessarily still acts on $\mathcal H^{\prime}$ by the same multiplier $m$. It will have a reproducing kernel $\mathcal K^\prime$  and our arguments starting with \eqref{eqn:4.1} remain valid for it, only the constant in \eqref{eqn:4.6} may be different. It follows that ${\mathcal K}^{\prime}=c^{\prime} \mathcal K$, hence  $\mathcal K_wv\in \mathcal H^{\prime}$ for all $w, v$, and therefore $\mathcal H^{\prime} =\mathcal H$.

\section{Two special cases and the invariant measure} In $\S\,3$ and $4$ we worked with a general representation $\pi$ of $K$ that was fixed throughout. Certain special choices of $\pi$ have important geometric interpretations. We first look at the case $\pi = \operatorname{Ad}_{\mathfrak p^+}$ which is an example of a unitary representation of $K$, with $V=\mathfrak p^{+}$. In this case we have
\begin{equation}\label{eqn:5.1}
\operatorname{Ad}_{\mathfrak p^+}(J(g, z))=g^{\prime}(z), 
\end{equation}
meaning the holomorphic tangent map (the complex Jacobian matrix) of $g$ acting on $D$. This is immediate for $g\in P^{+}$or $g \in K^{\mathbb C}$, as seen at the end of $\S\, 2$. To completely prove it, since both sides of \eqref{eqn:5.1} are multipliers, it is enough to also check it for $g \in P^{-}$ at $z=0$.
Now $g$ fixes $0$ and the holomorphic tangent space at $0$ can be identified with  $\mathfrak g^\mathbb C/ (\mathfrak k^{\mathbb{C}}+ \mathfrak p^{-})$. Writing $g=\exp Y$, $(Y \in \mathfrak p^{-})$ the tangent action is $\operatorname{ad}(Y)$, which is $I$ modulo $\mathfrak k^{\mathbb{C}}+ \mathfrak p^{-}$. This proves \eqref{eqn:5.1}. (A more direct proof can be found in \cite[II., Lemma 5.3]{S}.)

We observe that in this case (i.e., $\pi=\operatorname{Ad}_{\mathfrak p^+}$), \eqref{eqn:2.9}
becomes
$$
\operatorname{Ad}_{\mathfrak p^+} (J(a(t), 0))=
\exp \left(\frac{1}{2} \sum_1^r \log (1-x_j^2)\ad(h_j)\right)$$
which is diagonalized by the basis of $\mathfrak p^+$ formed by the root spaces for $\Delta_n^{+}$, the eigenvalues of $\ad(h_j)$ being $\gamma(h_j)$.

This is very useful for our second special case, namely the (scalar) representation $\pi=\operatorname{det} \operatorname{Ad}_{\mathfrak p^+}$. Now $\pi(J(g, z))=\operatorname{det}g^{\prime}(z)=\operatorname{Jac}_{g}(z)$, the complex Jacobian determinant. From \eqref{eqn:1.4} we immediately have
\begin{equation}\label{eqn:5.2} 
\operatorname{Jac}_{a(t)}(0)=\prod_j(1-x_j^2)^{p/2}.  
\end{equation}

Now we can get the ratio of the measures $\nu$ and $\lambda$. Since the real Jacobian is $|\operatorname{Jac}(z)|^2$, it is immediate from \eqref{eqn:3.4} that (with $Q(z)$ corresponding to our present $\pi$)
$Q(z)d\lambda(z)$ is invariant under the action of $G$. Therefore,
\begin{equation} \label{eqn:5.3}
d \nu(z)=c Q(z) d \lambda(z),
\end{equation}
and we also have, for $k \in K$, by (2.9) and by $\vert Jac_k(z)\vert ^2 = 1,$
\begin{equation} \label{eqn:5.4}
Q(ka(t), 0)=\prod_j (1-x_j^2)^p
\end{equation}

A few more things (not needed in $\S\,6$) follow easily. We can define the polynomial $h$ on $A\cdot 0$ by 
$$
h\left(\sum x_j e_j\right)=\prod_j\left(1-x_j^2\right)
$$
It is invariant under the Weyl group, so by a theorem of Chevalley (transferred
from $\mathfrak p$ to $\mathfrak p^+$) it uniquely extends to a $K$-invariant polynomial $h(z)$ on $\mathfrak p^+$. So we have 
\begin{equation} \label{eqn:5.5}
Q(z) = c h(z)^{-p}.
\end{equation}
Now consider the unitary representation $\pi^\prime = \pi^{-1}$ (with $\pi=\det \operatorname{Ad}_{\mathfrak p^+}$ as above). Corresponding to $\pi^\prime$, we have $Q^\prime=Q^{-1}$ and $\mathcal K^\prime = \mathcal K^{-1}$. In \eqref{eqn:3.2} working with $\pi^\prime$, and using that by\eqref{eqn:5.3} $Q^\prime$ and $Q$ cancel, the norm $\|f \|$ becomes the ordinary $L^2(D, \lambda)$-norm, and $\mathcal K^{\prime}(z, w)= \mathcal K(z, w)^{-1}$ becomes the classical Bergman kernel. As $\mathcal K$ was the polarization of $Q(z)^{-1}$, $\mathcal K^{\prime}$  is the polarization of $Q(z)$. Polarizing $h(z)$ gives a polynomial $h(z, w)$ in $z$ and $\overbar{w}$. By \eqref{eqn:5.5} therefore, the Bergman kernel of $D$ is
$$c \,h(z,w)^{-p}.$$
 as it was already proved in \cite[Sec. 3]{FK} and as we now proved it again (up to constant) in a rather indirect way. 
\section{The Harish-Chandra condition} We have seen at the end of $\S\, 4$  that $\mathcal H$, on which the discrete series representation acts, is non-zero if and only if
$$
\int_D\langle Q(z) v, v\rangle_V d \nu (z)< +\infty
$$
for all $v \in V$ ($Q$ is defined by \eqref{eqn:3.3}. We have to show that this is equivalent to the Harish-Chandra condition. First of all, our condition is equivalent to
\begin{equation} \label{eqn:6.1a}
\int_D Q(z) d \nu(z)<+\infty
\end{equation}
($\operatorname{GL}(V)$ - valued). By (3.5) and Schur's Lemma this is a scalar operator, but in order to keep the integrand symmetric, we write it in the equivalent form 
\begin{equation} \label{eqn:6.1}
\int_D \operatorname{tr}(Q(z)) d \nu(z) < +\infty.
\end{equation}
Now $\operatorname{tr}(Q(z))$ is a $K$-invariant function, so when we use the integral formula \eqref{eqn:2.10}, the $K$-integral disappears. Using also \eqref{eqn:5.3}, \eqref{eqn:5.4} to replace the measure $\nu$ with $\lambda$, i.e. with an ordinary integral in the $x$-space, we find that \eqref{eqn:6.1} is equivalent to
\begin{equation} \label{eqn:6.2}
\int_{0 \leq x_1 \leq \cdots  \leq x_r < 1} \operatorname{tr}(Q(x)) \prod_j \left(1-x_j^2\right)^{-p} P(x) d x_1 \cdots d x_r < +\infty. 
\end{equation}
Here $P(x)$ is the polynomial of $\S\, 2$, clearly irrelevant for the convergence of the integral.

To go on we have to describe the inducing representation $(\pi,V)$ more exactly. First of all, we denote by $\alpha_1, \ldots , \alpha_l$  the simple roots of $\mathfrak g^\mathbb C$, with $\alpha_1$ the only non-compact one (it is clear from $\mathfrak p^\mathbb{C} = \mathfrak p^+ \oplus \mathfrak p^-$ that there is only one.)

$\pi$ is a unitary representation of $K=K_{\text{ss}} \times Z$ (direct product, by simple connectedness). Therefore $\pi=\pi_0 \otimes \chi$ with $(\pi_0, V)$ irreducible unitary on $K_{\text{ss}}$, and $\chi$ a unitary character of $Z$. $\pi_0$ has a highest weight $\Lambda_0$ on $i\mathfrak h_{\text{ss}}$; we extend it to $i \mathfrak h$ by setting it to be $0$ on $h_{\alpha_1}$. To parametrize the characters of $Z$ in accordance with Harish-Chandra and \cite{W2}, we define the linear form $\Lambda_1$ on $i\mathfrak h$ by
$$
\Lambda_1\left(\ h_{\alpha_k}\right)= \begin{cases} 1 & (k=1) \\ 0 & (2 \leq k \leq l).\end{cases}
$$

Then, for $\lambda \in \mathbb{R}$, $\lambda \Lambda_1$ is the Lie algebra version of a character of $Z$. The linear function on $i\mathfrak h$
$$
\Lambda=\lambda_0 + \lambda \Lambda_1
$$
describes $\pi$ completely; it is  the highest weight of $\pi$.

The next step is to compute $\operatorname{tr}(Q(x))$ in \eqref{eqn:6.2}. We take a basis $\{v^s\}$ of $V$ consisting of weight vectors for $\pi_0$, with $v^{s}$ corresponding to the weight $\Lambda^s$.

From \eqref{eqn:2.9} we have, for each $v^s$,
\[\pi_0(J(a(t), 0)) v^s=\Big (\exp \frac{1}{2} \sum_j \log \left(1-x_j^2\right) \Lambda^s\left(h_j\right)\Big) v^s=\prod_j \left(1-x_j^2\right)^{\frac{1}{2} \Lambda^s\left(h_j\right)} v^s.\]
Similarly, with the character $\chi$ of $Z$ determined by $\lambda$, we have
$$
\chi(J(a(t), 0)) v^s=\prod_j(1-x_j^2)^{\frac{1}{2} \Lambda_1\left(h_j\right)} v^s. 
$$
It follows from the definition (3.3) of $Q$ that
$$
\operatorname{tr}(Q(x))=\sum_s \prod_j\left(1-x_i^2\right)^{-\left(\Lambda^s+\lambda \Lambda_1\right)\left(h_j\right)}.
$$
When we substitute this into \eqref{eqn:6.2} we get a sum (over $s$) of positive integrals. Their sum is finite if and only if each one is finite. It follows that the condition \eqref{eqn:6.2} is equivalent to
$$
-\left(\Lambda^s+\lambda \Lambda_1\right)\left(h_j\right) - p > - 1
$$
for all s and all $j$.  We have $\Lambda_1\left(h_j\right)=1$ for each $j$, by the definition of $\Lambda_1$ and using that $\alpha_1$ and each $\gamma_j$ have the same length (clear from looking at the two Cartan integers involving $\alpha_1$ and $\gamma_r$). So our condition can again be rewritten equivalently as
$$
\lambda<1-p-\Lambda^s\left(h_j\right) \quad(\forall s, \forall j).
$$

We will show that this condition is equivalent to the single inequality 
\begin{equation} \label{eqn:6.3}
\lambda<1-p-\Lambda_0\left(h_{r}\right).
\end{equation}
For this we have to prove that $\Lambda^s\left(h_j\right) \leq \Lambda_0\left(h_r\right)$, or what amounts to the same (since each $h_j$ has the same length),
\begin{equation} \label{eqn:6.4}
\left(\Lambda^s\mid \gamma_j\right) \leq \left(\Lambda_0 \mid \gamma_r\right)
\end{equation}
for all $s$ and $j$. We are in the real dual space of $i \mathfrak h=i \mathfrak z \oplus i \mathfrak h_{\text{ss}}$, and we prove (6.4) by using this
(orthogonal) decomposition. All $\Lambda^s$ are the same on $i \mathfrak z$, which is central in $\mathfrak k^\mathbb C$. Also each $\gamma_j$ is obviously the same on $i \mathfrak z$. So we only have to prove
\begin{equation} \label{eqn:6.5}
 (\overbar{\Lambda}^s\mid \overbar{\gamma}_j ) \leq \big (\overbar{\Lambda}_0 \mid \overbar{\gamma}_r\big )
\end{equation}
where the bars denote restriction to $i\mathfrak h_{\text{ss}}$. The $\bar{\gamma}_j$ are then the weights of the representation $\operatorname{Ad}_{\mathfrak p^{+}}$ of $\mathfrak{k}^{\mathbb C}_{\text{ss}}$; $\Lambda^s$ and $\Lambda_0$ were originally defined only on $i \mathfrak h_{\text{ss}}$ anyhow.

Observe that $\overbar{\Lambda}_0$ and $\overbar{\gamma}_r$ are the highest weights of two irreducible representations of $\mathfrak{k}_{ss}$ (namely $\pi_0$ and $\operatorname{Ad}_{\mathfrak p^{+}}$), while $\overbar{\Lambda}^{s}$ and $\overbar{\gamma}_{j}$ are other weights of the same two representations. It is a rather trivial general lemma that in such a situation an inequality like \eqref{eqn:6.5} holds. To prove it, we can find a Weyl group element $w$ such that $w\overbar{\Lambda}^s$ is dominant, hence equal to $\overbar{\Lambda}_0$. Then the thing to prove becomes $(\overbar{\Lambda}_0\mid w \overbar{\gamma}_j) \leq (\overbar{\Lambda}_0 \mid \overbar{\gamma}_r )$, which is trivial.

\eqref{eqn:6.3} is a form of the Harish-Chandra condition. In its original form \cite{HC1}, \cite{HC2} it was
$$
(\Lambda + \rho)(h_\gamma) < 0 \quad (\forall \gamma \in \Delta_n^{+}).
$$
It was shown in \cite[Theorem A.12]{RV}  by a direct elementary argument (easy to reconstruct using that every $\gamma$ is of the form $\gamma_r-\sum_2^{\ell} m_j \alpha_j$ with $m_j \geq 0$ and that $\rho(h_{\alpha_j})=1$) that it is equivalent to the single inequality
$$
(\Lambda+\rho) (h_r ) < 0.
$$
By \eqref{eqn:1.3} and by the definition of $\Lambda$ this is the same as \eqref{eqn:6.3}. \hfill $\Box$

\section{Final comments}
\begin{rem} \label{rem:1}The passage from $G_0$ to $G$ at the beginning of $\S\, 3$ is done by an \textit{ad hoc} method originally due to Satake. It restricts us to some extent and leaves open the question of the relation between $K$ (subgroup of $G$) and $K^\mathbb C$. Of course, the basic difficulty is that $G$ is not contained in any complex group that would be its complexification. There is, however, another way due to R. Herb and J, A. Wolf to deal with this situation \cite{HW}. They define a local group, which we tentatively denote by $G_{\text{loc}}^{\mathbb C}$, which contains $G$; putting it in place of $G^{\mathbb C}$ we can have our \eqref{eqn:2.1} and what follows from it right away in general with $G, K, Z, K^{\mathbb C}$ without having to use the subgroups of $G_0$ first and lift later.

It is done by taking $P^{+} K_0^{\mathbb C} P^{-} \subset G_0^{\mathbb{C}}$, and observing that $P^{+} K_0^{\mathbb C} P^{-}$ is a local group. Then one takes the universal cover $\xi: K^{\mathbb C} \rightarrow K_0$, the product $P^{+} \times K^{\mathbb C}\times P^{-}$and the map  $\operatorname{id} \times\operatorname{id} \xi \times \operatorname{id}$ mapping it to $P^{+} K_0^{\mathbb C} P^{-}$.  The product space will be $G_\text{loc}^{\mathbb C}$ after defining products on a sufficiently large subset of it. Naturally $P^\pm \times K^{\mathbb C}$ will have the products lifted from $P^\pm K_0^{\mathbb C}$ and will be subgroups of $G_0$. To define the product of an element of $P^{+}$ and an element of $P^{-}$ our 
identity \eqref{eqn:2.6} is useful: $\mathfrak p^-$ contains a mirror image $D^{-}$ of $D$, and products of elements in  $\exp D$ and in  $\exp D^{-}$ are defined via \eqref{eqn:2.6}. There is a little more work to do, it is all worked out carefully in \cite{HW}.
\end{rem}
\begin{rem}\label{rem:2} We defined the holomorphic discrete series, as customary in recent literature, by holomorphic induction from unitary representations $\pi$ of $K$. This is not the definition that Harish-Chandra uses \cite{HC1}, \cite{HC2}. He starts with a character $\chi$ of a torus $T$ (maximal in $K$ ) and does what amounts to holomorpphic induction directly to $G$.  This gives the same as our method, because of the theorem on "induction by stages". Holomorphic induction from $\chi$ to $K$ gives exactly all the irreducible unitary representations of $K$ (by the Borel- Weil theorem), and then inducing from $K$ to $G$ gives the same as inducing from $T$ directly to $G$.

Starting  with a character $\chi$ of $T$ has the advantage of immediately suggesting the results in \cite{HC2} about formal dimension and characters of the holomorphic discrete series. Such results are of course also present, and well known, for the representations of $K$; one could use them to give proofs of the results in \cite{HC2}.
The scalar given by our \eqref{eqn:6.1} is the formal dimension (cf. \cite[Lemma 2.10]{W2}) and suggests trying to directly evaluate such integrals.
\end{rem}

\end{document}